\documentclass[11pt]{amsart}
\usepackage[applemac]{inputenc}
\usepackage[american, frenchb]{babel}

\def\ladate{14 Août 2002}

\newif\ifarxiv
\arxivtrue

\newtheorem{theoreme}{Théorème}[]
\newtheorem{corollaire}[theoreme]{Corollaire}
\newtheorem{lemme}[theoreme]{Lemme}

\theoremstyle{definition}
\newtheorem{definition}{Definition}

\newtheorem{remarque}{Remarque}

\newcommand\NN{{\mathbb N}}

\newcommand\RR{{\mathbb R}}
\newcommand\ZZ{{\mathbb Z}}
\newcommand\HH{{\mathbb H}}

\newcommand\cA{{\mathcal A}}
\newcommand\cB{{\mathcal B}}
\newcommand\cE{{\mathcal E}}
\newcommand\cF{{\mathcal F}}

\newcommand\cS{{\mathcal S}}

\let\wh=\widehat
\let\wt=\widetilde

\renewcommand{\Re}{{\rm Re}}

\let\la=\lambda

\begin{document}

\ifarxiv\rightline{\ladate}\else
\textbf{Analyse Harmonique/\textit{Harmonic Analysis}}\hfill\ladate\break\fi
\ifarxiv\rightline{math/0208121}\fi

\bigskip

\title[Sur les Espaces de Sonine]{Sur les ``Espaces de Sonine''
associ\'es par de~Branges \`a la transformation de Fourier}

\author{Jean-Fran\c cois Burnol}

\ifarxiv\else\thanks{\textbf{Note présentée par J.-P. Kahane}}\fi

\begin{abstract}
Nous avons obtenu des formules explicites représentant les
fonctions $E(z)$ apparaissant dans la théorie des  ``Espaces de Sonine''
associés par De~Branges à la transformation de Fourier. 
\end{abstract}

\maketitle

\selectlanguage{american}
\title{\textbf{On the ``Sonine Spaces'' associated by de~Branges
to the Fourier Transform}}

\begin{abstract}
We have obtained explicit formulae representing the functions $E(z)$
arising in the theory of the ``Sonine Spaces'' associated by De~Branges
to the Fourier Transform.
\end{abstract}

\maketitle

\selectlanguage{french}
\bigskip

\begin{small}
Jean-Fran\c cois Burnol, Université Lille 1, UFR de Mathématiques, 
Cité Scientifique M2, F-59655 Villeneuve d'Ascq Cedex, France\\
\texttt{burnol@agat.univ-lille1.fr}
\end{small}

\parindent=0pt

\baselineskip = 16pt
\parskip=\the\baselineskip

%
\centerline{\hbox to4cm{\leaders \hrule height1pt \hfill}}
\selectlanguage{french}

\section{Espaces et Fonctions de de~Branges}

Nous commencerons par rappeler quelques éléments de base de la théorie
de de~Branges des ``Espaces de Hilbert de fonctions entières''
\cite{bra}. Tout ce dont nous aurons besoin est aussi exposé et démontré dans les
pages 220 à 228 du livre \cite{dymkean2} de Dym et McKean. Soit $L$
une droite dans le plan complexe, bordant un demi-plan ouvert $L_+$
qui est toujours chez de~Branges le demi-plan supérieur et qui sera
toujours chez nous le demi-plan $\Re(s) > 1/2$. Soit $w^\#$ le
symétrique orthogonal du nombre complexe $w$ par rapport à $L$, et
notons $F^\#(w) = \overline{F(w^\#)}$. Pour tout espace non nul de
Hilbert de fonctions entières tel que \textit{(a)} l'évaluation en
chaque $w$ est continue, \textit{(b)} $F\mapsto F^\#$ est une
isométrie (anti-unitaire), et \textit{(c)} si $F(\gamma) = 0$ alors
$G(w) = (w-\gamma^\#)/(w-\gamma) F(w)$ appartient à l'espace et est de
même norme que $F$, il existe par un théorème de de~Branges une (non
unique) fonction entière $E(w)$  telle que $|E(w)| > |E(w^\#)|$ pour
tout $w\in L_+$ et telle que l'espace de fonctions entières considéré
coïncide isométriquement avec l'espace $\cB(E)$ défini ainsi: $F$ est
dans $\cB(E)$ si $F(w)/E(w)$ est dans l'espace de Hardy $\HH^2(L_+)$
et si $F^\#(w)/E(w)$ l'est aussi. La norme hilbertienne est obtenu par
la première inclusion dans $\HH^2(L_+)$. Il faut préciser quel
multiple de la mesure de Lebesgue est pris sur la droite $L$: chez
de~Branges $L = \RR$ avec la mesure de Lebesgue, chez nous $L =
1/2 + i \RR$ avec $1/2\pi$ fois la mesure de Lebesgue. Un autre
théorème de de~Branges nous informe que $\cB(E)$ défini ainsi est
toujours non nul et est un espace de Hilbert, pour toute fonction
entière $E(w)$ satisfaisant à la condition ci-dessus. 
\let\ol = \overline

De plus on sait calculer exactement en fonction de $E$ les produits
scalaires entre évaluateurs. Dorénavant $L_+ = \{\Re(s)>1/2\}$. Alors
$w^\# = 1 - \ol{w}$. On supposera que la ``condition de réalité'' est
satisfaite, c'est-à-dire que pour tout $F(w)$ dans l'espace $\overline{F(\overline{w})}$
est dans l'espace avec la même norme.  On peut alors choisir
$E(w)$ de sorte que $\ol{E(w)} = E(\ol{w})$. On écrit $E = A - iB$
avec $A$ ``paire'' ($A(1-w) = A(w)$) et $B$ ``impaire''. Les fonctions
$A$ et $B$ ont tous leurs zéros sur la droite critique. Les autres
fonctions $E$ définissant le même espace (avec la même norme) et
réelles sur l'axe réel sont  connues exactement: il s'agit  des
fonctions $kA - iB/k$ pour $k\in \RR, k\neq 0$. 
La fonction (la condition de réalité est ici supposée vérifiée)
\begin{equation}
K(z_1, z_2) = \frac{E(z_1)E(z_2) - E(1-z_1)E(1-z_2)}{z_1 + z_2 -
1}
\end{equation}
\begin{equation}\label{noyau}
K(z_1, z_2) 
= 2\frac{(-iB(z_1))A(z_2)+A(z_1)(-iB(z_2))}{z_1 + z_2 - 1}
\end{equation}
est comme fonction de $z_2$ égale à l'évaluateur en $\ol{z_1}$ et
comme fonction de $z_1$ égale à l'évaluateur en $\ol{z_2}$.

\section{Espaces de Sonine}

\let\la=\lambda

Soit $\la>0$. Soit $\cF$ la transformation de Fourier
($\cF(\phi)(x) = \int_\RR \exp(2\pi ixy)\phi(y)dy$). Il est connu
classiquement (voir \cite[sect. 2.9]{dymkean1}) que $L^2(-\la,\la) +
\cF L^2(-\la,\la)$ est un sous espace \emph{fermé} de l'espace de
Hilbert $L^2(\RR)$ (cela se déduit aisément en particulier du
caractère compact de l'opérateur (à noyau) $F_\la = P_\la \cF P_\la$, où
$P_\la$ est la projection orthogonale sur $L^2(-\la,\la)$). Bien sûr
il s'agit d'un sous-espace propre et son complément perpendiculaire
est donc non-nul: il s'agit de l'espace de Hilbert des fonctions de
carrés intégrables, nulles et de Fourier nulles dans $(-\la,\la)$.
Suivant de~Branges \cite{bra}, nous dirons qu'une fonction est une
fonction de Sonine si elle satisfait à cette condition. 

L'argument le plus simple pour montrer concrètement l'existence de
fonctions de Sonine nous a été communiqué par le Professeur Kahane
\footnote{Lettre adressée à l'auteur, Mars 2002}: en appliquant un
polynôme approprié en $x$ et $d/dx$ à la distribution de Poisson
$\sum_{n\in\ZZ} \delta(x-n)$, on obtient une distribution tempérée
nulle et de Fourier nulle dans $(-A,A)$ avec $A$ arbitraire. Après
régularisation par multiplication et convolution additive par des
fonctions tests appropriées on obtient une fonction (non-nulle\dots)
dans la classe de Schwartz (donc de carré intégrable) ayant la
propriété de Sonine, avec $\la>0$ arbitrairement grand.

Dans \cite{jfcopoisson} nous montrons que la
\emph{convolution multiplicative} permet d'associer à toute distribution ayant
la propriété d'annulation de Sonine une formule de \emph{co-Poisson} portant sur
des fonctions de Sonine, et
qui généralise les formules de co-Poisson associées dans \cite{jfhab} à
la fonction dzêta de Riemann et aux séries $L$ de
Hecke-Tate. Il apparaît d'ailleurs fort utile pour l'étude des
fonctions de carrés intégrables de de~Branges-Sonine de disposer de
résultats généraux sur les distributions ayant la même propriété
d'annulation (propriété qu'il convient aussi de généraliser, par exemple en
remplaçant la condition ``nulle dans $(-\la,\la)$'' par la condition ``polynomiale dans
$(-\la,\la)$ de degré $\leq M$''). Dans la suite nous nous restreignons aux fonctions et
distributions \emph{paires}\footnote{Quelques ajustements simples partout dans ce qui suit
comme remplacer les cosinus par des sinus donnent la théorie pour le cas impair},
pour lesquelles $\cF$ devient la
transformée en cosinus $\cF_+(\phi)(x) = 2\int_0^\infty \cos(2\pi
xy)\phi(y)dy$. Nous désignons par $K_\la$ l'espace des fonctions de
Sonine paires et de carrés intégrables. Nous utiliserons le résultat
suivant:

\begin{theoreme}[{\cite{jfcopoisson}}]
Soit $D$ une distribution tempérée paire qui s'annule sur
$(-\la,\la)$. On (rappelle que l'on) peut donner un sens à $\wh{D}(s) = \int_0^\infty
D(t)t^{-s}dt$ comme fonction analytique de $s$ pour $\Re(s)\gg 1$. Si
$\cF_+(D)$ est à nouveau nulle dans $(-\la,\la)$ alors $\wh{D}(s)$ est
une fonction entière de $s$, qui a des zéros triviaux en $-2n$,
$n\in\NN$. De plus on a l'équation fonctionnelle 
$$\pi^{-s/2}\Gamma(s/2)\wh{D}(s) =
\pi^{-(1-s)/2}\Gamma((1-s)/2)\wh{\cF_+(D)}(1-s)$$
\end{theoreme}

Considérons alors l'espace de Hilbert non-nul $\cS_\la$ des fonctions
entières $F(w) = \pi^{-w/2}\Gamma(w/2)\wh{f}(w)$, où $f$ parcourt
$K_\la$. Il est aisé de vérifier que $\cS_\la$ vérifie les axiomes
\textit{(a)}, \textit{(b)}, et \textit{(c)}, et est donc un espace de
de~Branges, par rapport à la droite critique, et qui vérifie la
condition de réalité (De~Branges \cite{bra64, bra} a associé plus
généralement de tels espaces aux transformations de Hankel de
paramètre $\nu$, $\nu>-1$). Le problème (qui est équivalent à celui du
calcul des produits scalaires entre évaluateurs dans $K_\la$) était donc posé de déterminer
les fonctions entières $\cA_\la(w)$, $\cB_\la(w)$ et $\cE_\la(w) =
\cA_\la(w) - i \cB_\la(w)$ qui permettent d'écrire l'identité
(y-compris pour la norme)  $\cS_\la = \cB(\cE_\la)$. Nous donnons dans
la présente Note une solution explicite à ce problème.

\section{Un problème de projection orthogonale}

Soit
$P_\la$ la projection orthogonale sur $L^2(-\la,\la)$ et $\wt{P_\la}=
\cF P_\la \cF^* = \cF^* P_\la \cF$ la projection sur
$\cF(L^2(-\la,\la))$. Soit $L^2(\RR)^{\rm pair}$ l'espace des fonctions paires 
avec comme norme $\|f\|^2=\int_0^\infty |f(t)|^2dt$. 
Sur $L^2(\RR)^{\rm pair}$ on a  $\wt{P_\la} =
\cF_+ P_\la \cF_+$. Soit 
$K_\la$ son sous-espace de Sonine, et $G_\la =
P_\la(L^2(\RR)^{\rm pair}) + \wt{P_\la}(L^2(\RR)^{\rm pair})$. 

\begin{lemme}
L'espace de Sonine $K_\la$ est le noyau de l'opérateur $P_\la + \wt{P_\la}$ (sur 
$L^2(\RR)^{\rm pair}$). La 
restriction de cet opérateur à $G_\la$ est 
auto-adjoint et de Fredholm (i.e. de la forme 1+compact) et est inversible. La
projection orthogonale de $f\in L^2(\RR)^{\rm pair}$ sur $G_\la$ est
$(P_\la + \wt{P_\la})^{-1} (P_\la(f) + \wt{P_\la}(f))$.
\end{lemme}

Si $P_\la(f) + \wt{P_\la}(f) = 0$ alors $g=P_\la(f)=-\wt{P_\la}(f) = 0$
 et $f$ est dans $K_\la$.
Donc si $k\in G_\la$ a la même image sous $P_\la + \wt{P_\la}$ que $f$
alors $f-k$ est dans $K_\la$ et $k$ est la projection orthogonale de $f$ sur $G_\la$.
Considérons $\phi: L^2(-\la,\la)^{\rm pair}\oplus
L^2(-\la,\la)^{\rm pair} \to G_\la$ qui envoie $(u,v)$ sur $u +
\cF_+(v)$. L'opérateur borné $\phi$ est injectif et surjectif donc un
isomorphisme topologique qui conjugue $P_\la + \wt{P_\la}$ à
l'opérateur (lui-aussi auto-adjoint) $(u,v) \mapsto (u + F_\la v,
F_\la u + v)$, avec $F_\la = P_\la \cF_+ P_\la$ compact auto-adjoint
de norme strictement inférieure à $1$. L'inverse de cet opérateur est
$(u,v) \mapsto ((1 - D_\la)^{-1} (u - F_\la v), (1 - D_\la)^{-1} (-
F_\la u + v))$. Dans cette formule  $D_\la = F_\la^2$  peut-être
remplacé par le noyau $P_\la \cF^* P_\la \cF P_\la$ puisqu'il n'agit que sur
des fonctions paires. On note que $F_\la$ a les mêmes vecteurs propres
que $D_\la$: ce sont les ``fonctions prolates sphéroidales'' (paires)
$e_{2n}$, $n\geq0$ (\cite{slepian}). Ainsi:

\begin{lemme} 
Les vecteurs propres de l'opérateur (dans $G_\la$) 
de Fredholm $P_\la
+ \wt{P_\la}$ sont les vecteurs $e_{2n} \pm \cF_+(e_{2n})$, $n\in\NN$. 
Ils forment une base orthogonale de $G_\la$ puisque l'opérateur est
auto-adjoint.
\end{lemme}

Il est donc possible d'écrire la projection orthogonale sur l'espace
de Sonine en utilisant la base orthogonale de $G_\la$ ci-dessus, mais
nous ne ferons pas usage de telles formules.

\begin{remarque} 
L'auteur doit à M.~Balazard et É.~Saias (communication privée,
décembre 2001) l'observation que les fonctions sphéroidales per\-mettent
d'écrire une base orthogonale de $G_\la$ (leur base orthogonale
diffère légèrement de celle apparaissant ici).
\end{remarque}

Si l'on combine les formules précédentes on obtient:

\begin{theoreme}
Soit $f(x)$ une fonction paire de carré intégrable. Sa projection
orthogonale sur l'espace de Sonine $K_\la$ est donnée par la formule:
$$\pi_\la(f) = f - (1 - D_\la)^{-1}\Big(P_\la(f) - F_\la\cF_+(f)\Big)
- \cF_+(1 - D_\la)^{-1}\Big(P_\la\cF_+(f) -  F_\la(f)\Big)$$
On a $F_\la = P_\la \cF_+ P_\la$ et $D_\la = F_\la^2$ 
agissant sur $L^2(-\la,\la)^{\rm pair}$. On peut remplacer $D_\la$ par le noyau 
$\sin(2\pi\la(x-y))/\pi(x-y)$ restreint à $L^2(-\la,\la)^{\rm pair}$.
\end{theoreme}

\begin{corollaire}
Si $f|_{(-\la,\la)} = 0$ alors 
$$\pi_\la(f) = f\Big|_{|t|>\la} - \Big(\cF_+(1 -
D_\la)^{-1}P_\la\cF_+(f)\Big)\Big|_{|t|>\la}$$
Si $f = \cF_+(f)$ alors
$$\pi_\la(f) = f - (1 + \cF_+)(1+F_\la)^{-1}P_\la(f)$$
Si $f = -\cF_+(f)$ alors
$$\pi_\la(f) = f - (1 - \cF_+)(1-F_\la)^{-1}P_\la(f)$$
\end{corollaire}

\section{Deux remarquables distributions ayant la propriété de Sonine}

Nous appliquons formellement la projection orthogonale aux
distributions paires invariante et anti-invariante sous Fourier:
$\delta_{-\la}(t) + \delta_\la(t) + 2\cos(2\pi\la t)$ et
$\delta_{-\la}(t) + \delta_\la(t) - 2\cos(2\pi\la t)$ et cela nous
mène à 

\begin{definition}
On pose
$$A_\la(t) = \delta_{-\la}(t) + \delta_\la(t) + 2\cos(2\pi\la t) -
\Big((1 + \cF_+)(1+F_\la)^{-1}(2\cos(2\pi\la y))\Big)(t)$$
$$-iB_\la(t) = \delta_{-\la}(t) + \delta_\la(t) - 2\cos(2\pi\la t) +
\Big((1 - \cF_+)(1-F_\la)^{-1}(2\cos(2\pi\la y))\Big)(t)$$
\end{definition}

\begin{theoreme}
Les distributions paires tempérées $A_\la(t)$ et $-iB_\la(t)$ ont la
propriété de Sonine: elles sont nulles et de Fourier nulles dans
l'intervalle $(-\la,\la)$. Elles sont respectivement invariante et
anti-invariante sous Fourier.
\end{theoreme}

\begin{definition}
On pose
$$\psi_+^\la(t) = 2\cos(2\pi\la t) -
\cF_+\Big((1+F_\la)^{-1}P_\la(2\cos(2\pi\la y))\Big)(t)$$
$$\psi_-^\la(t) = 2\cos(2\pi\la t) +
\cF_+\Big((1-F_\la)^{-1}P_\la(2\cos(2\pi\la y))\Big)(t)$$
\end{definition}

\begin{theoreme}
La fonction entière paire $\psi_\pm^\la$ est (l'unique) fonction
entière dont la restriction à $(-\la,\la)$ est $(1\pm
F_\la)^{-1}(2\cos(2\pi\la y))$. Elle est aussi \break{l'unique} fonction
solution de l'équation $\phi(x) \pm \cF(\phi|_{(-\la,\la)})(x) =
2\cos(2\pi \la x)$. On a $A_\la = \psi_+^\la + \cF(\psi_+^\la)$ et
$-iB_\la = - \psi_-^\la + \cF(\psi_-^\la)$. 
\end{theoreme}

\section{Les fonctions $\cE(w)$, $\cA(w)$, $\cB(w)$ pour les espaces
de Sonine}

Nous avons démontré le résultat suivant:

\begin{theoreme}
Soit 
$$\cE_\la(w) = \pi^{-w/2}\Gamma(w/2)\Big(\la^{1/2 - w} +
\frac{\sqrt{\la}}{2}\int_\la^\infty (\psi_+^\la(t) -
\psi_-^\la(t))t^{-w}dt\Big)$$
L'intégrale est absolument convergente pour $\Re(w)>0$. La fonction
$\cE_\la(w)$ est une fonction entière satisfaisant la condition de
de~Branges pour le demi-plan $\Re(s)>1/2$ et l'espace $\cB(\cE_\la)$
coïncide (isométriquement) avec l'espace $\cS_\la$ des transformées de
Mellin complétées des fonctions de $K_\la$.
\end{theoreme}

Nous esquissons la démonstration. On notera $X_w^\la(t)$ l'unique
élément de $K_\la$ tel que $\forall f\in K_\la\ \int_0^\infty
f(t)t^{-w}dt = \int_0^\infty f(t)X_w(t)dt$ (pour $\Re(w)\leq 1/2$ la
première intégrale est un prolongement analytique). Pour $\Re(w)>1/2$
la formule pour la projection orthogonale permet d'écrire
``explicitement'' $X_w^\la(t)$, et on voit en particulier que
$X_w^\la(t)- t^{-w}$ est la restriction à $t>\la$ d'une fonction
entière. L'idée est de calculer la distribution $(td/dt +
w)X_w^\la(t)$: en effet on peut montrer en partant de la formule
\ref{noyau} que (la transformée de Mellin complétée de) la partie
invariante sous Fourier de cette distribution sera un multiple (réel
positif lorsque $w\in(1/2,+\infty)$) d'une fonction $\cA_\la(z)$
possible, et que (la transformée de Mellin complétée de) la partie
anti-invariante sous Fourier sera un multiple d'une fonction
$-i\cB_\la(z)$.  Or le calcul montre que la partie invariante de
$(td/dt + w)X_w^\la(t)$ est multiple de la distribution remarquable
$A_\la(t)$ définie précédemment et que la partie anti-invariante est
multiple de $-iB_\la(t)$. Une étude plus poussée permet d'affirmer que
la combinaison $\sqrt{\la}/2 (A_\la - iB_\la)(t)$ a comme transformée
de Mellin (complétée par le facteur Gamma) une fonction $\cE_\la(w)$
possible pour l'espace de Hilbert de de~Branges Sonine $\cS_\la$. La
démonstration donne aussi:

\begin{theoreme}
Il y a coïncidence entre $\cE_\la(w)$ et $\sqrt{\la}$ fois la valeur
du ``saut'' de l'évaluateur (``euclidien'') $Z_w^\la(t) =
\pi^{-w/2}\Gamma(w/2)X_w^\la(t)$ en $t= \la$.
\end{theoreme}

Il est intéressant que certaines quantités apparaissant dans cette
étude et liées au noyau de Dirichlet jouent un rôle dans l'étude de
son déterminant de Fredholm sur un intervalle fini, déterminant dont
on connaît l'importance dans la théorie des matrices aléatoires
\cite{mehta}.


\bibliographystyle{amsplain}

\end{document}